\documentclass[11pt]{article}
\usepackage{tikz}
\usepackage{subfigure}
\usepackage[english]{babel}
\usetikzlibrary{arrows, graphs}
\usepackage[center]{caption2}
\usepackage{amsfonts,amssymb,amsmath,latexsym,amsthm}
\usepackage{multirow}
\usepackage{bbm}
\usepackage{thmtools}
\usepackage{thm-restate}%定理重述
\usepackage[colorlinks=true,allcolors=blue,backref=page]{hyperref}
\usepackage[noabbrev,capitalize,nameinlink]{cleveref}

\newcommand \CC{\mathcal C}
\newcommand \HH{\mathcal H}
\newcommand \sdeg{\operatorname{sdeg}}
\newtheorem{thm}{Theorem}[section]
\newtheorem{cor}[thm]{Corollary}
\newtheorem{lem}[thm]{Lemma}

\newtheorem{conj}[thm]{Conjecture}

\newtheorem{prob}[thm]{Problem}
\newtheorem{obs}[thm]{Observation}
\newtheorem{claim}{Claim}

\begin{document}

\title{Unavoidable induced subgraphs forced by graphs with many vertices of prescribed properties}
\author{Jin Sun$^a$,\quad Xinmin Hou$^{b,c}$\footnote{E-mail: jinsun@mail.ustc.edu.cn (J. Sun),  xmhou@ustc.edu.cn (X. Hou)}\\
\small $^a$School of Mathematical Sciences\\
\small Anhui University, Hefei, Anhui 230621, China\\
\small$^b$ School of Mathematical Sciences,\\
\small University of Science and Technology of China, Hefei 230026, Anhui, China\\
\small$^c$ Hefei National Laboratory,\\
\small University of Science and Technology of China, Hefei 230088, Anhui, China\\
}

\date{}

\maketitle 	
\begin{abstract}
Given a function $p : V(G)\to \mathbb N$ and an integer $k\ge 0$, define $p_k(G)$ as the number of vertices with $p(v)\ge k$. 
We say that $p_k(G)$ is bounded for all $\HH$-free graphs if there exists a constant $c=c(\HH)$ such that $p_k(G)<c$ for all such graphs $G$. Here, a graph $G$ is said to be $\HH$-free if it contains no member of $\HH$ as an induced subgraph.
When $p$ represents the degree of a vertex, Ramsey's theorem implies that $p_0(G)$ is bounded for all $\{K_n, \overline K_n\}$-free graphs, where $K_n$ and $\overline K_n$ denote the complete graph and the edgeless graph on $n$ vertices, respectively.   
The connected version of Ramsey's theorem says that $p_0(G)$ is bounded for all $\{K_n, P_n, K_{1,n}\}$-free connected graphs, where $P_n$ and $K_{1,n}$ are the $n$-vertex path and the star with $n$ leaves.
In this paper, we extend Ramsey's theorem to $p_2(G)$ where $p$ denotes the degree, the local independence number, the local component number, and sharp degree, that is,  we characterize the forbidden family of graphs $\HH$ such that $p_2(G)$ is bounded for all (connected) $\HH$-free graphs. Moreover, we also characterize the forbidden family of graphs $\HH$ for which there is a constant $c=c(\HH)$ such that $p_c(G)$ is bounded for all $\HH$-free graphs.  
\end{abstract}

\section{Introduction}
For positive integers $m$ and $n$ with $m\le n$, let $[m,n]$ denote the set $\{m,m+1,\dots ,n-1,n\}$, and write $[n]$ for $[1,n]$.
In this paper, we consider only finite, undirected and simple graphs. Let $G=(V(G),E(G))$ be a graph. For vertices $y,z\in V(G)$, we denote the adjacency of $y$ and $z$ by $y\sim z$.   Let  $\mathbbm 1_{y\sim z}$ be the {\it characteristic function} which equals $1$ if $y$ is adjacent to $z$, $0$ otherwise.
For a vertex $v\in V(G)$, let $N_G(v)=\{u\in V(G): u\sim v\}$ denote the {\it neighborhood} of $v$ in $G$ (the subscript may be omitted if there is no confusion).
For a subset $S\subseteq V(G)$, we denote by $G[S]$ the subgraph of $G$ induced by $S$. Let $\alpha(S)$ be the {\it stability number} (or {\it independence number}) of the induced subgraph $G[S]$, and $c(S)$ be the {\it number of connected components} of $G[S]$. For two disjoint subsets $A, B$ of $V(G)$, we denote by $G[A,B]$ the bipartite subgraph induced by edges between $A$ and $B$. For additional notation, see \cite{d}.

For graphs $H_1$ and $H_2$, we say $H_1\prec H_2$ if $H_2$ contains an induced subgraph isomorphic to $H_1$. For a family $\HH$ of graphs, we say $G$ is {\it $\HH$-free}, if there is no graph $H\in \HH$ such that $H\prec G$. %The family  $\HH$ is called forbidden subgraphs.
For families $\HH_{1}$ and $\HH_2$, we say $\HH_1\le \HH_2 $ if for any graph $H_2\in\HH_2$, there exists an $H_1\in \HH_1$ such that $H_1\prec H_2$. 
Throughout this article, when we say a graph $G$ contains an $H$, or we have an $H$ in $G$, we mean that $G$ contains an induced subgraph isomorphic to $H$, i.e., $H\prec G$. 
The following straightforward result will be commonly used without explicit mention.
\begin{obs}\label{folklore}
The relation `$\le$' is transitive. If $\HH_1\le \HH_2$, then every $\HH_1$-free graph is also $\HH_2$-free. 
\end{obs} 

In 1929, the well-known Ramsey's theorem appeared.

\begin{thm}(Ramsey's Theorem, \cite r) 
For any positive integers $m$ and $n$, there exists a minimum positive integer $R=R_{m}(n)$ such that if the edges of the complete graph $K_R$ are colored with $m$ colors, then there is a monochromatic clique of order $n$, i.e., a clique all of whose edges have the same color.
\end{thm}

Specifically, for $m=2$, Ramsey's theorem says that for any large $n$, a graph with $R_{2}(n)$ vertices must contain  $K_n$ or $\overline K_n$ as an induced subgraph, where $\overline K_n$ is the edgeless graph of order $n$, i.e., the complementary graph of $K_n$. 

Given a graph class $\mathcal C$, define 
%\[B(\CC)=\{\HH:\text{there is a constant } c=c(\HH)  \text{ such that } |V(G)|<c \text{ for all } \HH \text{-free graphs } G\in \CC \}. \]
\[B(\CC)=\{\HH: \text{ there exists $c=c(\HH)$ such that } |V(G)|<c \text{ for all } \HH \text{-free graphs } G\in \CC \}. \]
A Ramsey-type problem is to determine $B(\CC)$ when $\CC$ is a class of graphs with additional property. 
The following is a connected version of the classical Ramsey's theorem.
	\begin{thm}[Prop 9.4.1 in \cite{d}]\label{THM:Conn-Ramsey}
There is a minimum integer $N_0(n)$ such that every connected $\{K_n\,,K_{1,n}\,,P_n\}$-free  graph $G$ satisfies $|V(G)|<N_0(n)$.   
  \end{thm}
Let $\CC$ be the family of connected graphs. Thus according to \cref{THM:Conn-Ramsey}, $\HH\in B(\CC)$ if and only if $\HH\le \{K_n,K_{1,n},P_n\}$ for some positive integer $n$ . 
We list some results for different graph families in the following. (We do not write down explicitly the forbidden subgraph conditions, as there are too many members to list. The interested readers are encouraged to refer to the corresponding papers.)

\begin{description}

	\item[(1)] Allred, Ding, and Oporowski \cite{ado} determined $B(\CC)$,where $\CC$ is the family of 2-connected graphs.
	\item[(2)] Ding and Chen \cite {dc} determined $B(\CC)$, where $\CC$ is the family of doubly connected graphs, i.e., connected graphs whose complement graphs are also connected.
	\item[(3)] Chudnovsky, Kim, Oum, and Seymour \cite{ckos} determined $B(\CC)$, where $\CC$ is the family of prime graphs. 
\end{description}

Given a graph parameter $\mu$, define 
\[B(\mu)=\{\HH: \text{there is a constant } c=c(\HH) \text{ such that } \mu(G)<c \text{ for all } \HH \text{-free graphs } G \}\]
It is very interesting to determine $B(\mu)$ for different graph parameters $\mu$. We list several findings  in the literature.
\begin{description}
	%\item[(1)] Furuya \cite{f} determined $B(\gamma_c)$, where $\gamma_c$ is the connected domination number of a graph. %if and only if $\HH \le \{K_n^*,K_{1,n}^*,P_n\}$ for some positive integer $n$, where $\gamma_c$ is the connected domination number of a graph. In other words, there exists a minimum integer $\gamma_c(n)$ such that each connected $\{K_n^*,K_{1,n}^*,P_n\}$-free graph $G$ satisfies $\gamma_C(G)<\gamma_c(n)$.

	%\item[(2)] Chiba and Furuya \cite{cf} determined $B(\mu)$ when $\mu$ is the path cover/partition number.
	
	%\item[(3)] Choi, Furuya, Kim and Park \cite{cfkp} determined $B(\nu)$, where $\nu(G)$ is the matching number of $G$.

	\item[(1)] Galvin, Rival, and Sands \cite{grs}, and Atminas, Lozin, and Razgon \cite{alr} determined $B(\mu)$, where $\mu(G)$ is the length of the longest path of graph $G$.
    \item[(2)] Kierstead and Penrice \cite{kp}, and Scott, Seymour, and Spirkl \cite{sss} determined $B(\delta)$, where $\delta(G)$ is the minimum degree of graph $G$.
	\item[(3)] Lozin and Razgon \cite{lr} determined $B(\text{tw})$, where $\text{tw}(G)$ is the tree-width of graph $G$.
    \item[(4)] Sun and Hou \cite{sh} determined $B(\text{def})$, where $\text{def}(G)$ is the deficiency of graph $G$.
	%\item[(3)] Lozin \cite{l} determined $B(\mu)$, where $\mu$ is the neighborhood diversity or VC-dimension.
	
   %\HH必须是有限集
	\end{description}
Since there exist infinitely many graphs with sufficiently large chromatic number and
sufficiently large girth~\cite{e}, the following conjecture in this flavor is equivalent to the Gy\'arf\'as-Sumner conjecture~\cite{g,s81}.
\begin{conj}%[{\bf Gy\'arf\'as-Sumner Conjecture}]
Let $\HH$ be a finite family of connected graphs. Then $\HH\in B(\chi)$ if and only if $\HH\le \{K_n,T\}$ for some tree $T$ and positive integer $n$, where $\chi(G)$ is the chromatic number of graph $G$. 
%For given positive integer $n$ and tree $T$, there exists a constant $c=c(n,T)$ such that all $\{K_n,T\}$-free graphs satisfy $\chi (G)<c$, where $\chi(G)$ is the chromatic number of graph $G$.    
\end{conj}

Given a function $p : V(G)\to \mathbb N$ and an integer $k\ge 0$, define 
$$p_k(G)=\#\{v\in V(G): p(v)\ge k\}.$$
For example, $p(v)$ can be $\deg(v)$, $\alpha(N(v))$, and $c(N(v))$ etc. 
We say that $p_k(G)$ is bounded for $\HH$-free graphs $G$ if there exists a constant $c=c(\HH)$ such that $p_k(G)<c$ for $\HH$-free graphs $G$.

Define 
\begin{align*}
&B_1(p)=\{\HH:  p_2(G) \text{ is bounded  for all connected } \HH \text{-free graphs } G\}, \\[1mm]
&B_2(p)=\{\HH: p_2(G) \text{ is bounded for all } \HH \text{-free graphs } G\}, \text{ and}\\[1mm]
&B_3(p)=\{\HH: \text{ there is $c=c(\HH)$ such that $p_{c}(G)$ is bounded for all } \HH \text{-free graphs } G\}.
\end{align*}
Therefore, a graph family $\HH\in B_3(p)$ means there are 
constants $c_1=c_1(\HH)$ and  $c_2=c_2(\HH)$ such that $p_{c_1}(G)<c_2$ for all $\HH$-free graphs $G$. For convenience, we may assume $c_1,c_2\ge 2$ throughout the paper.

\begin{prob}\label{main prob}
Given a function $p:V(G)\to \mathbb N$, determine $B_i(p)$ for $i=1,2,3$.	
\end{prob}

%Specifically, $p(v)$ will be $\deg(v)$, 
Let $G$ be a graph and a vertex  $v\in V(G)$.  Define the {\it local independence number} of $v$ as $\alpha_L(v)=\alpha(N(v))$, and 
the {\it local component  number} of $v$ as $c_L(v)=c(N(v))$.
The {\it sharp degree} of $v$ is defined as $\sdeg(v):=c(G-v)-c(G)+1$. Thus $v$ is a {\it cut vertex} if and only if $\operatorname{sdeg}(v)\ge 2$. Therefore, the sharp degree of a vertex can be viewed as a generalization of cut vertex. Note that for a vertex $v\in V(G)$ , we have $\deg(v)\ge \alpha_L(v)\ge c_L(v)\ge \sdeg(v)$. 
In this article,  we determine $B_i(p)$ for $i=1,2,3$, where $p$ represents $\deg$, $\alpha_L$, $c_L$, and $\sdeg$, respectively.

In order to state our Ramsey-type results, we introduce some kinds of graphs. As usual, we let $K_n, \overline K_n, P_n$ denote the {\it complete graph}, {\it edgeless graph} and {\it path} of order $n$, respectively, and let $K_{s,t}$ be the {\it complete bipartite graph} with partitions of orders $s$ and $t$. For two graphs $G_1$ and $G_2$, we define the {\it join} $G_1+G_2$ by the graph with vertex set $V(G_1)\cup V(G_2) $ and edge set $E(G_1)\cup E(G_2)\cup \{ xy :x\in V(G_1),\, y\in V(G_2) \} $. We denote by $nG$ the graph consisting of $n$ disjoint copies of $G$. The graph $G^+$ obtained by adding a pendant vertex to vertex $v$ of graph $G$ is formally defined as $V(G^+)=V(G)\cup \{v'\}$ and $E(G^+)=E(G)\cup \{vv'\}$. 
The following are more kinds of graphs (See \cref{fig 1}).  
\begin{itemize}
	\item $K_{1,n}^{*} $: The graph obtained by adding a pendant vertex to each leaf of $K_{1,n}$.
	\item $K_n^{*}$: The graph obtained by adding a pendant vertex to each vertex of $K_n$.
	\item $K_n \times K_2$: The graph obtained by adding a perfect matching between two disjoint copies of $K_n$.
	\item $G_n$: The graph obtained by adding $n$ pendant vertices to each vertex of $K_n$.
\end{itemize}

\begin{figure}
\centering
\begin{tikzpicture}[baseline=10pt,scale=0.85]
\filldraw[fill=black,draw=black] (0,0) circle (0.1);
\filldraw[fill=black,draw=black] (1,2) circle (0.1);
\filldraw[fill=black,draw=black] (1,1) circle (0.1);
\node at (1,-0.2) {$\vdots$};
\filldraw[fill=black,draw=black] (1,-1.5) circle (0.1);
\filldraw[fill=black,draw=black] (2,2) circle (0.1);
\filldraw[fill=black,draw=black] (2,1) circle (0.1);
\filldraw[fill=black,draw=black] (2,-1.5) circle (0.1);
\draw (0,0) -- (1,2);
\draw (0,0) -- (1,1);
\draw (0,0) -- (1,-1.5);
\draw (2,2) -- (1,2);
\draw (2,1) -- (1,1);
\draw (2,-1.5) -- (1,-1.5);
\node at (1,-3) {$K_{1,n}^*$};

\draw (3.5,-1.8) rectangle (4.5,2.3);
\filldraw[fill=black,draw=black] (4,-1.5) circle (0.1);
\filldraw[fill=black,draw=black] (4,1) circle (0.1);
\filldraw[fill=black,draw=black] (4,2) circle (0.1);
\node at (5,-0.2) {$\vdots$};
\node at (4,-0.3) {$K_n$};
\filldraw[fill=black,draw=black] (5,-1.5) circle (0.1);
\filldraw[fill=black,draw=black] (5,1) circle (0.1);
\filldraw[fill=black,draw=black] (5,2) circle (0.1);
\draw (4,-1.5) -- (5,-1.5);
\draw (4,1) -- (5,1);
\draw (4,2) -- (5,2);
\node at (4,-3) {$K_{n}^*$};

\draw (7.5,-1.8) rectangle (6.5,2.3);
\draw (8.5,-1.8) rectangle (9.5,2.3);
\filldraw[fill=black,draw=black] (7,-1.5) circle (0.1);
\filldraw[fill=black,draw=black] (9,-1.5) circle (0.1);
\filldraw[fill=black,draw=black] (7,1) circle (0.1);
\filldraw[fill=black,draw=black] (9,1) circle (0.1);
\filldraw[fill=black,draw=black] (7,2) circle (0.1);
\filldraw[fill=black,draw=black] (9,2) circle (0.1);
\draw (7,-1.5) -- (9,-1.5);
\draw (7,1) -- (9,1);
\draw (7,2) -- (9,2);
\node at (9,-0.2) {$\vdots$};
\node at (7,-0.3) {$K_n$};
\node at (8,-3) {$K_n \times K_2$};

\draw (13,0) circle (1.5);
\node at (13,0) {$K_n$};
\filldraw(13,1) circle (0.1);
\filldraw(12.5,2) circle (0.1);
\filldraw(13,2) circle (0.1);
\node at (13.5,2) {$\cdots$};
\filldraw(14,2) circle (0.1);
\draw (13,1) --  (12.5,2);
\draw (13,1) --  (13,2);
\draw (13,1)--   (14,2);
\node at (11,0) {$\vdots$};
\node at (15,0) {$\vdots$};
\filldraw(13,-1) circle (0.1);
\filldraw(12.5,-2) circle (0.1);
\filldraw(13,-2) circle (0.1);
\node at (13.5,-2) {$\cdots$};
\filldraw(14,-2) circle (0.1);
\draw (13,-1) --  (12.5,-2);
\draw (13,-1) --  (13,-2);
\draw (13,-1)--   (14,-2);
\node at (13,-3) {$G_n$};
\end{tikzpicture}
\caption{The graphs $K_{1,n}^*,\, K_n^*,\, K_n \times K_2$ and $G_n.$}
\label{fig 1}
\end{figure}

Now, we are ready to state our results for $B_1(p)$
%when $p$ is the degree, the local stable number, the local component number, and  $\deg, $
on \cref{main prob}.

\begin{restatable}{thm}{csdeg}
\label{deg} 
$\HH\in B_1(\deg)$ if and only if \[\HH \le \{K_n\,, P_n\,, K_{1,n}^*\,, K_{2,n}\,, K_2+nK_1\,, K_1+nK_2 \} \]
 for some positive integer $n$.
\end{restatable}

\begin{restatable}{thm}{localstablenumber} \label{local stable number}
$\HH\in B_1(\alpha_L)$
if and only if 
$$ \HH\le \{K_n^*\,, P_n\,, K_{1,n}^*\,, K_{2,n}\,, \overline K_2+K_n\,, K_1+nP_3\,, K_n \times K_2 \} $$ for some positive integer $n$.
\end{restatable}

\begin{restatable}{thm}{localconnected} \label{local connected}
$\HH \in B_1(c_L)$ if and only if 
 $ \HH\le \{K_n^*\,, K_{1,n}^*\,, P_n\,, K_{2,n}\,, K_n \times K_2\} $
 for some positive integer $n$.
\end{restatable}

\begin{restatable}{thm}{cutvertex} \label{cut vertex}
$\HH\in B_1(\sdeg)$
if and only if $\HH\le \{K_n^*\,, K_{1,n}^*\,, P_n\} $ 
for some positive integer $n$.
\end{restatable}

As corresponding corollaries of Theorems \ref{deg}, \ref{local stable number}, \ref{local connected} and \ref{cut vertex}, the following are our results for $B_2(p)$ on \cref{main prob}.

\begin{cor} 
$\HH\in B_2(\deg)$ if and only if 
\[\HH\le \{K_n, nP_3, nK_3, K_{1,n}^*\,, K_{2,n}\,, K_2+nK_1, K_1+nK_2 \} \] for some positive integer $n$.
\end{cor}

\begin{proof}
We first prove the ``only if'' part. Let $\HH\in B_2(\deg)$. Thus there is a constant $c=c(\HH)$ such that $(\deg)_2(G)< c $ for every  $\HH$-free graph $G$. 
%The number of vertices with degree at least $2$ for 
Clearly,
$$ (\deg)_2(nP_3)=n,\,(\deg)_2(K_{1,n}^*)=n+1,$$ and $$(\deg)_2(H)=|V(H)| \text{ for } H\in\{K_n, nK_3, K_{2,n}, K_2+nK_1, K_1+nK_2\}.$$
Let $n=c$. Then the number of vertices with degree at least $2$ for all these graphs is not bounded by $c$, implying that they are not $\HH$-free. Hence $\HH\le \{K_n, nP_3, nK_3, K_{1,n}^*\,, K_{2,n}\,, K_2+nK_1, K_1+nK_2 \} $.

Next we prove the ``if'' part. 
Assume
$$\HH\le \{K_n, nP_3, nK_3, K_{1,n}^*\,, K_{2,n}\,, K_2+nK_1, K_1+nK_2 \}. $$ 
By \cref{folklore}, every $\HH$-free graph $G$ is also $\{K_n, nP_3, nK_3, K_{1,n}^*\,, K_{2,n}\,, K_2+nK_1, K_1+nK_2 \} $-free. 
Let $m$ be the number of connected components of $G$ which contain a vertex of degree at least $2$. Denote these $m$ connected components by $\{C_i:i\in [m]\}$, where $u_i\in C_i$ and $v_i,w_i\in N(u_i)$. Thus for each $i\in [m]$, $C_i$ contains an induced $K_3$ or $P_3$ according to whether $v_i$ is adjacent to $w_i$.
Since $G$ is $\{nP_3, nK_3\}$-free, $m\le 2n-2$. Applying \cref{deg} to the connected $\{K_{4n}, P_{4n}, K_{1,4n}^*\,, K_{2,4n}\,, K_2+4nK_1, K_1+4nK_2 \} $-free graph $G[C_i]$ for each $i\in [m]$, there exists a constant $c=c(n)$, which depends only on $n$, such that  
 $(\deg)_2(G[C_i])<c$.
Hence \[(\deg)_2(G)=\sum_{i=1}^m (\deg)_2(G[C_i])<(2n-2)c.\]

\end{proof}

By a similar technique, we can get the following. 
\begin{cor}
$\HH\in B_2(\alpha_L)$
if and only if 
$$ \HH\le \{K_n^*, nP_3, K_{1,n}^*\,, K_{2,n}\,, \overline K_2+K_n, K_n \times K_2 \} $$ for some positive integer $n$.
\end{cor}

\begin{cor}
$\HH\in B_2(c_L)$
if and only if 
$$ \HH\le \{K_n^*, nP_3, K_{1,n}^*\,, K_{2,n}\,, K_n \times K_2\} $$for some positive integer $n$.
\end{cor}

\begin{cor}
$\HH\in B_2(\sdeg)$
if and only if 
 $ \HH\le \{K_n^*\,, nP_3, K_{1,n}^*\} $ for some positive integer $n$.
\end{cor}

We characterize  $B_3(p)$ where $p$ is $\deg$, $\alpha_L$, $c_L$, and $\sdeg$, respectively, in the following.

\begin{restatable}{thm}{finitedeg} \label{finite deg}
$\HH\in B_3(\deg)$
if and only if $\HH\le \{K_n, K_{n,n}\,, nK_{1,n}\}$ for some positive integer $n$.
\end{restatable}

\begin{restatable}{thm}{finitelocalstablenumber} \label{finite local stable number}
$\HH\in B_3(\alpha_L)$ if and only if $\HH \le \{K_{n,n}\,, nK_{1,n}\,, K_n+\overline K_n\,, G_n\}$ for some positive integer $n$.
\end{restatable}

\begin{restatable}{thm}{finitelocalcomponentnumber}\label{finite local component number}
$\HH\in B_3(c_L)$ if and only if $\HH \le \{K_{n,n}\,, nK_{1,n}\,, G_n\}$ for some positive integer $n$.	
\end{restatable}

\begin{restatable}{thm}{finitesdeg} \label{finite sdeg}
$\HH \in B_3(\sdeg)$
if and only if $\HH \le \{nK_{1,n}, G_n\}$ for some positive integer $n$.
\end{restatable}
For a function $p:V(G)\to \mathbb N$, the {\it H-index}~\cite{es,hhbgm} $H_p$ of $p$ is the maximum integer $h$ such that $p_h(G)\ge h$. Using the corresponding H-index $H_p$ of $p$, one can translate $B_3(p)$ into $B(H_p)$.

%For a scholar with $n$ publications $\{v_i : i\in [n]\}$, let $p(v_i)$ be the number of citations to the publication $v_i$. Thus the H-index of $p$ is the scholar's H-index to evaluate the academic influence of this scholar. We have the following observation.

\begin{obs}\label{zhuanhua}
$B_3(p)=B(H_p)$, where $H_p$ is the H-index of $p$.	
\end{obs}

\begin{proof}
On the one hand, let $\HH\in B_3(p)$. Thus there exist constants $c_1=c_1(\HH)$ and $c_2=c_2(\HH)$  such that  $p_{c_1}(G)<c_2$  for any  $\HH$-free graph $G$. Let $c=c_1+c_2$. Thus 
$$p_{c}(G)\le p_{c_1}(G)<c_2<c,$$
implying that $H_p(G)<c$. Therefore, $\HH\in B(H_p)$.

On the other hand, let $\HH\in B(H_p)$. Thus there exists a constant $c=c(\HH)$ such that $H_p(G)<c$  for any  $\HH$-free graph $G$, i.e., $p_c(G)<c$. Let $c_1=c_2=c$, then $p_{c_1}(G)<c_2$. Hence $\HH\in B_3(p)$.	

Therefore, $B_3(p)=B(H_p)$.
\end{proof}

According to \cref{zhuanhua}, Theorems \ref{finite deg}, \ref{finite local stable number}, \ref{finite local component number} and \ref{finite sdeg} can be viewed as determining $B(H_p)$, where $H_p$ is the corresponding H-index of $p$. 

The rest of this paper is arranged as follows. In section $2$, we will determine $B_1(p)$ and prove Theorems \ref{deg}, \ref{local stable number}, \ref{local connected} and \ref{cut vertex}. 
In section $3$, we will determine $B_3(p)$ and prove Theorems \ref{finite deg}, \ref{finite local stable number}, \ref{finite local component number} and \ref{finite sdeg}.

\section{Proofs of Theorems \ref{deg}, \ref{local stable number}, \ref{local connected} and \ref{cut vertex}}
%An {\it induced matching} $M$ of a graph $G$ is an induced 1‐regular subgraph of $G$. The induced matching number of $G$ is defined as $$\nu'(G) = \max\{|M| : M \text{ is an induced matching of $G$}\}.$$

Using the notation $\prec$ defined earlier, the {\it induced matching number} of $G$ is defined as 
$$\nu'(G) = \max\{n : nK_2\prec G\}.$$
The following technical lemma has been proved by the authors in~\cite{sh}.
\begin{lem}[Lemma 2.4 in~\cite{sh}]\label{induced matching}
Let $G$ be a bipartite graph with vertex partition sets $X$ and $Y$. Suppose $\deg(x)\ge 1$ and $\deg(y)\le n$ for all $x\in X,\, y\in Y$.
%$\delta(X)\ge 1$ and $\Delta(Y)\le n$.
If  $|X|\ge n(p-1)+1 $, then $\nu' (G)\ge p$.
%, where $\nu'(G)$ is the induced matching number of $G$. 
Moreover, the lower bound of $|X|$ is tight.
\end{lem}

%\begin{proof}
%$X\subseteq N(Y)$ as $\deg(x)\ge 1\; \forall x\in X$. Take minimal subset $Y'$ of $Y$ dominates $X$, i.e., $N(Y')=X$ and for any proper subset $Y''$ of $Y'$, $N(Y'')\not =X$. So $Y'$ is irredundant, i.e., every $y\in Y'$ has a private neighbor $x_y\in X$ such that $N(x_y)\cap Y'=\{y\}$. Since each vertex in $Y'$ has at most $n$ neighbors in $X$ and $|X|\ge n(p-1)+1$, $|Y'|\ge p$. Choose $p$ vertices in $Y'$, and a private neighbor for each one, they form an induced $pK_2$. So $\nu' (G)\ge p$.  

%Specially, let $G=(X,Y,E)$ be a bipartite graph with $X=\{x_{ij}:i\in [p-1],\, j\in [n]\}$, $Y=\{y_i:i \in [p-1]\}$, $E=\{(x_{ij},y_i):i \in [p-1],\, j\in [n]\}$. It satisfies $\deg(x)\ge 1,\, \deg(y)\le n,\; \forall x\in X\,, y\in Y$ but $\nu' (G)= p-1$ since $|X|$ only equals to $ n(p-1) $.
%\end{proof}

We restate \cref{deg} as follows.
\csdeg*
\begin{proof}%[Proof of  Theorem \ref{deg}]{}
 
We first prove the ``only if'' part. Assume that $\HH\in B_1(\deg)$, i.e., there is a constant $c=c(\HH)$ such that 
$(\deg)_2(G)< c$ for every connected $\HH$-free graph $G$.  
Clearly, \begin{align*}
(\deg)_2(K_n)&=n, &(\deg)_2(K_{2,n})=n+2,\qquad (\deg)_2(K_2+nK_1)&=n+2,\\ 
(\deg)_2(P_{n})&=n-2, &(\deg)_2(K_{1,n}^*)=n+1, \qquad (\deg)_2(K_1+nK_2)&=2n+1.
\end{align*}
 Let $n= c+3$, then all these graphs are connected but the number of vertices with degree at least $2$ is not bounded by $c$. Hence they are not $\HH$-free. Therefore, $\HH\le \{K_n, P_n, K_{1,n}^*, K_{2,n}, K_2+nK_1, K_1+nK_2 \} $.

Next we prove the ``if'' part. Now suppose $\HH\le \{K_n, P_n, K_{1,n}^*, K_{2,n}, K_2+nK_1, K_1+nK_2 \} $. Since an $\HH$-free graph is also $\{K_n, P_n, K_{1,n}^*, K_{2,n}, K_2+nK_1, K_1+nK_2 \}$-free,  it suffices to show that if a connected graph $G$ satisfies 
$$(\deg)_2(G)\ge N_0(N_1), \text{ where }N_1=(n-1)R_2(2n-1),$$   
then it contains one of $K_n, P_n, K_{1,n}^*, K_{2,n}, K_2+nK_1$ or  $K_1+nK_2$ as an induced subgraph. Here $N_0(N_1)$ is the constant depending only on $N_1$  from Theorem~\ref{THM:Conn-Ramsey}. Suppose $G$ is such a graph. For a cut vertex $v$ in an induced subgraph $H$ in $G$,
it has at least $2$ neighbors in $H$, and thus $\deg_G(v)\ge 2$. Deleting all the vertices of degree $1$ of $G$ all at once, 
the induced subgraph $G'$ left with order at least $(\deg)_2(G)\ge N_0(N_1)$ is still connected. By Theorem~\ref{THM:Conn-Ramsey}, $G'$ contains $K_{N_1}, P_{N_{1}}$ or $K_{1,N_1}$ as an induced subgraph. As $N_1\ge n$, we are done if $G'$ contains $K_{N_1}$, or $P_{N_{1}}$.

Now assume $G'$ contains a copy of $K=K_{1,N_1}$. Denote the central vertex of $K$ by $s$, and the leaves by $X=\{x_i: i\in [N_1]\}$. Since $\deg_G(x_i)\ge 2$ in $G$, we can find another vertex $y_i\not =s$ in $G$ adjacent to $x_i$. (Note that $y_i$ and $y_j$, $i,j\in [N_1],\, i\not =j$ may be the same vertex.) If there is a $y_i$, satisfying $|N(y_i)\cap X|\ge n$, then $\{s,y_i\}\cup (N(y_i)\cap X)$ contains a copy of $K_{2,n}$ or a copy of $K_2+nK_1$ as an induced subgraph according to whether $y_i$ is adjacent to $s$.

Therefore, we may assume $|N(y_i)\cap X|<n$ for all $i\in [N_1].$ Let $Y=\{y_i: i\in[N_1]\}$. Consider the bipartite graph $B=G[X,Y]$. We have $|X|=N_1=(n-1)R_2(2n-1)$, $\deg_{B}(x)\ge 1$, and $\deg_B(y)\le n-1$ for all $x\in X$ and $y\in Y$. Applying \cref{induced matching} to the bipartite graph $B$, we can find $X'\subseteq X,\, Y'\subseteq Y$ which form an $R_2(2n-1)K_2$ in $G[X,Y]$. If $Y'$ contains a clique of order $n$, then $G$ contains a $K_n$. We are done. Otherwise $Y'$ must contain a stable set $S$ of order $2n-1$ as $|Y'|=R_2(2n-1)$. If $|N(s)\cap S|\ge n$, then $G$ contains an $K_1+nK_2$ as an induced subgraph. Otherwise, by the pigeonhole principle, $|N(s)\cap S|< n$, there are at least $n$ vertices not adjacent to $s$, thus $G$ contains $K_{1,n}^*$.

In other words, every connected $\{K_n, P_n, K_{1,n}^*, K_{2,n}, K_2+nK_1, K_1+nK_2 \} $-free graph $G$ satisfies $(\deg)_2(G)\le c=N_0(N_1)$. 
The proof is completed. 
 \end{proof}

\cref{local stable number} is restated as follows.
\localstablenumber*

\begin{proof}%[Proof of Theorem \ref{local stable number}]
For a graph $G$, the vertex $v$ such that $\alpha(N(v))=1$ is not a cut vertex of any induced subgraph of $G$ which contains $v$. Deleting all these vertices all at once, the induced subgraph $G'$ left is still connected. 
 
We first prove the ``only if'' part. Assume that $\HH\in B_1(\alpha_L)$, i.e., there is a constant $c=c(\HH)$ such that 
$(\alpha_L)_2(G)< c$ for every connected $\HH$-free graph $G$. When $G=K_n^*, P_n, K_{1,n}^*, K_{2,n}, \overline K_2+K_n, K_1+nP_3$, or  $K_n \times K_2$, the order of $G'$ is  $n, n-2, n+1, n+2, n, n+1, 2n$ respectively. Let $n= c+3$. Then $(\alpha_L)_2(G)=|V(G')|> c$. Hence all these graphs $G$ are not $\HH$-free. Thus $ \HH\le \{K_n^*, P_n, K_{1,n}^*, K_{2,n}, \overline K_2+K_n, K_1+nP_3, K_n \times K_2 \} $.

Next we prove the ``if'' part. Set 
\[N_3=R_{2^8}(n+2),\  N_2=nR_2(n)N_3, \text{ and } N_1=2N_2-1. \]
All we need to show is that if a connected graph $G$ satisfies 
$$(\alpha_L)_2(G)\ge N_0(N_1),$$
then it contains at least one of $\{K_n^*, P_n, K_{1,n}^*, K_{2,n}, \overline K_2+K_n, K_1+nP_3,K_n \times K_2\}$ as an induced subgraph.
Since the connected graph $G'$ satisfies $|V(G')|= (\alpha_L)_2(G)\ge N_0(N_1)$, according to  Theorem~\ref{THM:Conn-Ramsey}, $G'$ contains $P_{N_1},\, K_{N_1}$ or $K_{1,N_1}$ as an induced subgraph. 
If $G'$ contains $P_{N_1}$, then we are done as $P_n\prec P_{N_1}$.
We divide the following proof into two cases. 

\vspace{5pt}
\noindent{\bf Case 1:} $G'$ contains a clique $X$ of order $N_1$.

Note that all vertices in $K$ have local stability number at least $2$ in $G$. For  a vertex $x\in X$, if no vertex $y\in V(G)$ is adjacent to $x$ but not complete to  $X$, then to ensure $\alpha_L(x)\ge 2$, there must exist non-adjacent vertices  $y$ and $z$ that are complete to $X$. 
Hence $G$ contains $\overline K_2+K_n$.

Now assume that for any $x\in X$, there exists a $y$ adjacent to $x$ but not complete to $X$. Denote by $Y$ the set consisting of these corresponding $y$.
%$X=\{x_i: i\in [N_1]\}$, and $\forall x_i, \exists y_i$ adjacent to $x_i$ but not complete to $X$. 
If there is a vertex $y\in Y$ such that $|N(y)\cap X|\ge n$, then $G$ contains an $\overline K_2+K_n$, formed by the vertex $y$, a vertex $z\in X$ not adjacent to $y$, and $n$ vertices in $N(y)\cap X$, done.
Thus we may assume  $|N(y)\cap X|< n$ for all $y\in Y$.  Applying Lemma \ref{induced matching} to the bipartite graph $G[X,Y]$, there exists an induced matching $\{x_iy_i: i\in [R_2(n)]\}$ in $G[X,Y]$, where $x_i\in X$ and  $y_i\in Y$.   If $Y'=\{y_i: i\in [R_2(n)]\}$ contains a clique of order $n$, then $G$ contains a $K_n \times K_2$; otherwise, $Y'$ contains a stable set of order $n$, which implies a $K_n^*$ in $G$.

\vspace{5pt}

\noindent{\bf Case 2:} $G'$ contains an induced $K_{1,N_1}$. 

Denote the central vertex by $s$, and  the set of leaves by $X$. For $x\in X$, the condition $\alpha_L(x)\ge 2$ is ensured if either  there exists a neighbor $y$ of $x$ not adjacent to $s$ (we denote the set of all such $x$ by $X_1$), 
or $N(x)\subseteq \{s\}\cup N(s)$ and there exists non-adjacent vertices $y$ and $z$ in $N(x)\cap N(s)=N(x)-\{s\}$. (we denote the set of all such $x$ by $X_2$). Consequently, $X$ is the disjoint union of $X_1$ and $X_2$.
As $N_1=2N_2-1$, we have two subcases by the pigeonhole principle. 

\noindent{\bf Subcase 2.1:} %At least $N_2$ vertices in $X$ satisfies  it has a neighbor not adjacent to $s$. 
 $|X_1|\ge N_2$.

Denote by $Y$ the set consisting of neighbors of $x$ that are not adjacent to $s$ for  $x\in X_1$. If there exists a $y\in Y$ satisfying $|N(y)\cap X_1|\ge n$, then $G$ contains $K_{2,n}$ (formed by those vertices $y, s$ and $n$ vertices in  $N(y)\cap X_1$). Otherwise,  applying Lemma \ref{induced matching} to the  bipartite graph $G[X_1,Y]$, there exists an induced matching $\{x_iy_i: i\in [R_2(n)]\}$ in $G[X_1,Y]$, where $x_i\in X_1$ and $y_i\in Y$.
If $Y'=\{y_i:i \in [R_2(n)]\}$ contains a clique of order $n$, then we have a  $K_n^*$ in $G$; otherwise, $Y'$ contains a stable set of order $n$, then we have a $K_{1,n}^*$ in $G$.

\noindent{\bf Subcase 2.2:} $|X_2|\ge N_2$.

Hence we can find vertices $\{x_i, y_i, z_i: i\in [N_2]\}$ such that 
$\{s, x_i, y_i, z_i\}$ forms an induced $K_4-e$, where all $x_i\in X_2$ are distinct and $y_iz_i\notin E(G)$. Note that $y_i$ and $z_i$ are distinct and symmetric. However, $y_i$ may be the same vertex as $y_j$ or $z_j$ for $i\not=j$. If some vertex $s'$ indeed appears $nR_2(n)$ times, we replace $s$ with $s'$ and return to the Subcase 2.1. For example, if $s'=z_i$ for all $i\in [nR_2(n)]$, then $\{s', x_i, y_i: i\in [nR_2(n)]\}$ is the desired structure for Subcase 2.1.

Therefore, without loss of generality, we may assume all vertices $y_i, z_i$ for $i\in [N_3]$ are distinct. 
Consider the auxiliary complete graph $K_{N_3}$. Color the edge $(i,j),\, i<j$ by  the $\{0,1\}$-sequence
$$(\mathbbm 1_{x_i\sim y_j},\,
\mathbbm 1_{x_i\sim z_j},\,
\mathbbm 1_{y_i\sim x_j},\,
\mathbbm 1_{z_i\sim x_j},\,
\mathbbm 1_{y_i\sim y_j},\,
\mathbbm 1_{z_i\sim z_j},\,
\mathbbm 1_{y_i\sim z_j},\, 
\mathbbm 1_{z_i\sim y_j} ).$$
Note that $N_3=R_{2^8}(n+2)$. By Ramsey's theorem, there exists a monochromatic clique $K$ of order $n+2$ in the $2^8$-edge-coloring of $K_{N_3}$.  Without loss of generality, we may assume $V(K)=[n+2]$. Now we consider the color of $K$.

For $i,j \in [n+2],\, i<j$, suppose $\mathbbm 1_{y_i\sim y_j}=1$ first. If $\mathbbm 1_{x_i\sim y_j}=1$, then $\{x_1, x_2, y_i: i\in [3,n+2]\}$ forms an induced $\overline K_2+K_n$. Otherwise, by symmetry  ($y_i$ and $z_i$ are symmetric,  $i$ and $j$ are symmetric, e.g. the arguments for $\mathbbm 1_{x_i\sim y_j}$ and $\mathbbm 1_{y_i\sim x_j}$ are symmetric), we have  $\mathbbm 1_{x_i\sim y_j}= \mathbbm 1_{y_i\sim x_j}=0$. Then $\{x_i, y_i: i\in [n]\}$ forms an induced $K_n^*$.
Now assume $\mathbbm 1_{y_i\sim y_j}=\mathbbm 1_{z_i\sim z_j}=0$ for $ i<j$ with $i,j \in [n+2]$. If $\mathbbm 1_{x_i\sim y_j}=1$, then  $\{x_1, x_2, y_i: i\in [3,n+2]\}$ forms an induced $K_{2,n}$. If $\mathbbm 1_{y_i\sim z_j}=1$, then  $\{y_1, y_2, z_i: i\in [3,n+2]\}$ forms an induced $K_{2,n}$. 
Therefore, by symmetry, we can assume 
\[\mathbbm 1_{x_i\sim y_j}=\mathbbm 1_{y_i\sim x_j}=\mathbbm 1_{x_i\sim z_j}=\mathbbm 1_{z_i\sim x_j}=\mathbbm 1_{y_i\sim z_j}=\mathbbm 1_{z_i\sim y_j}=0.\] 
Hence $\{s, x_i, y_i, z_i: i\in [n]\}$ forms an induced $K_1+nP_3$.
\end{proof}

Now we give the proof of Theorem~\ref{local connected}. Furuya \cite{f} characterized the forbidden induced subgraph family for the {\it connected domination number}. Let $\gamma_c(G)$ be the connected domination number of graph $G$, i.e., the minimum size of a connected subset $D\subseteq V(G)$ satisfying $D\cup \bigcup_{d\in D} N(d)=V(G)$. 
\begin{thm}\cite{f} \label{domination number}
For every positive integer $n$, there exists a minimum integer $\gamma_c(n)$ such that every connected $\{K_n^*\,, K_{1,n}^*\,,P_n\}$-free graph $G$ satisfies $\gamma_c(G)<\gamma_c(n)$.
\end{thm}
We remark that although the conclusion in Furuya \cite{f} is stated in terms of the domination number, the proof actually constructs a connected dominating set of size less than $\gamma_c(n)$ for every connected $\{K_n^*\,, K_{1,n}^*\,,P_n\}$-free graph.

We restate \cref{local connected} as follows.
\localconnected*

\begin{proof}  %[Proof of Theorem \ref{local connected}]
We first prove the ``only if'' part. Assume $\HH\in B_1(c_L)$. Thus there exists a constant $c=c(\HH)$ such that $(c_L)_2(G)<c$ for all connected $\HH$-free graph $G$.
Given a graph $G$, set the corresponding $G'=G[\{v\in V(G):c_L(v) \ge 2\}]$.
For $G=K_n^*\,, P_n\,, K_{1,n}^*\,, K_{2,n}$ or $K_n \times K_2$, the corresponding $G'$ has order $n, n-2, n+1, n+2, 2n$ respectively. Let $n=c+2$. Then $(c_L)_2(G)= |G'|\ge c$. This implies that
$$ \HH\le \{K_n^*\,, P_n\,, K_{1,n}^*\,, K_{2,n}\,, K_n \times K_2\}.$$

Next we prove the ``if'' part. All we need to show is that if a connected graph $G$ satisfies 
$$(c_L)_2(G)>  \gamma_c(n)R_2(N_1), \text{ where }  N_1=nN_2,\, N_2=R_2(n),$$
then it contains $K_n^*\,, P_n\,, K_{1,n}^*\,, K_{2,n}$, or  $K_n \times K_2$ as an induced subgraph. 

Let $D$ be a minimum dominating set of $G$. If $|D|\ge \gamma_c(n)$, then by \cref{domination number}, $G$ contains at least one of $K_n^*, P_n$, or $K_{1,n}^*$ as an induced subgraph. We are done. Thus we may assume $|D|< \gamma_c(n)$. Let $G'=G[\{v:c_L(v)\ge 2\}]$. Thus $|G'|=(c_L)_2(G)> \gamma_c(n)R_2(N_1)$. Therefore,  
\[|V(G')|-|D|> \gamma_c(n)(R_2(N_1)-1).\]
By the pigeonhole principle,  there exists a vertex $d\in D$ such that $|N(d)\cap V(G')|\ge R_2(N_1)$. 
If $N(d)\cap V(G')$ contains a clique of order $N_1$, then we have a $K_{N_1}$ in $G'$. Otherwise, $N(d)\cap V(G')$ contains a stable set of order $N_1$, thus we find a star $K_{1,N_1}$ with all leaves in $G'$.

\noindent{\bf Case 1:} $G'$ contains $K=K_{N_1}$ with $V(K)=\{x_i: i\in [N_1]\}$. 

Since every vertex $v\in V(G')$ satisfies $c_L(v)\ge 2$ in $G$, there exists a vertex $y_i\in N(x_i)$ such that $y_i$ cannot be connected to $V(K)-\{x_i\}$ in $G-\{x_i\}$ for each $x_i\in V(K)$. Moreover, $y_i\not =y_j$ if $i\not =j$, for all $i,j\in [N_1]$. 
Since $N_1\ge R_2(n)$, by Ramsey's Theorem,  the subgraph induced by $Y=\{y_i: i\in [N_1]\}$  contains a clique of order $n$ or a stable set of order $n$. The former case implies that $G$ contains a $K_n \times K_2$, and the latter case implies that $G$ contains a $K_n^*$.

\noindent{\bf Case 2:} There is a $K=K_{1,N_1}$ with $V(K)=\{s, \,x_i: i\in [N_1]\}$ such that $s$ is the central vertex, and $x_i\in V(G')$, for all $i\in [N_1]$. 

Since $c_L(v)\ge 2$ in $G$ for all vertices $v\in V(G')$, there exists a neighbor $y_i$ of $x_i$ but not adjacent to $s$ for all $i\in [N_1]$. If for some $i$, $|N(y_i)\cap K|\ge n$, then $G[\{s,y_i,\}\cup (N(y_i)\cap K)]$ contains a $K_{2,n}$, done. Otherwise, by Lemma \ref{induced matching}, we may assume $N(y_i)\cap K=\{x_i\}$ for $i\in [N_2]$.  If $Y=\{y_i: i\in [N_2]\}$ contains a $K_n$, then $G$ contains a $K_n^*$; otherwise, $Y$ contains an $\overline K_n$, and hence $G$ contains a $K_{1,n}^*$.
\end{proof}

We restate \cref{cut vertex} as follows.
\cutvertex*

\begin{proof}
We first prove the ``only if'' part. 
Recall that $(\sdeg)_2(G)$ is the number of cut vertices of $G$. Hence $(\sdeg)_2(K_n^*)=n, (\sdeg)_2(P_n)=n-2$, and $(\sdeg)_2(K_{1,n}^*)=n+1$. Hence $\HH\le \{K_n^*\,, K_{1,n}^*\,, P_n\}$ for $n=c+2$.

Next we prove the ``if'' part. By Theorem \ref{domination number}, there is a constant $\gamma_c(n)$ such that every connected  $\{K_n^*, K_{1,n}^*, P_n\}$-free graph $G$ satisfies $\gamma_c(G)\le \gamma_c(n)$.

\begin{claim}\label{CL:cut} 
Let $D$ be a connected dominating set of connected graph $G$, then $D$ contains all cut vertices of $G$.
%For a connected graph $G$, the connected dominating set $D$ contains all cut vertices of $G$.
\end{claim}
\begin{proof}Suppose not. Then there exists a cut vertex $v$ of $G$ such that $v\not \in D$. Hence $G-v$ is not connected. Since $D\subseteq V(G-v)$ and $G[D]$ is connected, there exist a connected component $C$ of $G-v$ that is disjoint with $D$. This contradicts the fact that $D$ is a dominating set.
\end{proof}
By Claim~\ref{CL:cut}, the number of cut vertices of connected $\HH$-free graph $G$ is bounded by $\gamma_c(n)$.
%So if $\HH\le \{K_n^*, K_{1,n}^*, P_n\} $, then every connected $\HH$-free graph $G$ is also $\{K_n^*, K_{1,n}^*, P_n\}$-free, the number of cut vertices of $G$ is bounded by $\gamma_c(n)$. 

\end{proof}

\section{Proofs of Theorem \ref{finite deg}, \ref{finite local stable number}, \ref{finite local component number} and \ref{finite sdeg}}
%The following lemma is similar to Ramsey's theorem. As for proof, we may assume $k=2$, and consider the bipartite case. Then use the same technique as the classical proof of Ramsey's theorem.
As folklore, the following Ramsey-type lemma is needed.
\begin{lem}\cite{alr} \label{bipartite Ramesy }
 For any positive integer $k$ and $q$, there exists a number $M\!R(k,q)$ such that in every $k$-partite graph $G=(V_{1},V_{2},\dots , V_{k}, E) $ with $|V_i|\ge M\!R(k,q)$ for all $i\in [k]$, there is a collection of subsets $U_i\subseteq V_i$ of order $|U_i|=q$ satisfying every pair of subsets  induces either a $K_{q,q}$ or $E_{2q}$.
 \end{lem}

We restate \cref{finite deg} as follows.
\finitedeg*
\begin{proof}
We first prove the ``only if'' part. When $n=c_1+c_2$, it holds that $(\deg)_{c_1}(K_n)=n$, $(\deg)_{c_1}(K_{n,n})=2n$, and $(\deg)_{c_1}(nK_{1,n})=n$.
%, the number of vertices with degree at least $c_1$ is $n, 2n, n$ respectively. 
Thus they are not $\HH$-free, which implies that  $\HH\le \{K_n, K_{n,n}, nK_{1,n}\}$.

Next we prove the ``if'' part. 

\begin{claim}Let $N_3=R_2(n)$, $N_2=R_{2^{n^2+2n+1}}(2n)$, and $N_1=N_2\cdot N_3+N_2$.  For every $\{K_n, K_{n,n}, nK_{1,n}\}$-free graph $G$, it holds that
$$(\deg)_{N_1}(G)<N_2.$$
%\[\#\{v\in V(G):\deg(v)\ge N_1\} <N_2. \]
\end{claim} 

\begin{proof}
Suppose not. Let  $\{v_i: i\in [N_2]\}$ be a set of vertices with degrees at least $N_1$. As $N_1=N_2\cdot N_3+N_2$, for each $v_i$, we can find vertices $\{v_i^j: j\in [N_3]\}$ within $ N(v_i)$. Moreover, all $v_i^j$ with $i\in [N_2],\, j\in [N_3]$  are distinct even though $v_i$ is adjacent to $v_{i'}$ and $v_{i'}^j$ as $\deg(v_i)\ge N_2\cdot N_3+N_2$. For every $i$, we can find $n$ stable vertices from $\{v_i^j: j\in [N_3]\}$ since there is no $K_n$. For convenience, we denote these $n$ stable vertices by $\{v_i^j: j\in [n]\}$.  
Color the auxiliary complete graph $K_{N_2}$ by $2^{n^2+2n+1}$ colors as follows. The edge $(i,i'),\, i<i'$ is colored by 
\[ (\mathbbm 1_{v_i\sim v_{i'}}, \mathbbm 1_{v_i\sim v_{i'}^1}, \cdots, \mathbbm 1_{v_i\sim v_{i'}^n}, \mathbbm 1_{v_{i'}\sim v_{i}^1},\cdots, \mathbbm 1_{v_i'\sim v_{i}^n}, \mathbbm 1_{v_i^1\sim v_{i'}^{1}}, \cdots, \mathbbm 1_{v_i^j\sim v_{i'}^{j'}},\cdots ,\mathbbm 1_{v_i^n\sim v_{i'}^n}),\]
where  $j,j'\in [n]$.
By Ramsey's Theorem, there exists a monochromatic clique $K$ of order $2n$ in the auxiliary complete graph $K_{N_2}$ as $N_2=R_{2^{n^2+2n+1}}(2n)$. For convenience, suppose $V(K)=[2n]$. Consider the color of $K$.

Since $G$ is $K_n$-free, we must have $\mathbbm 1_{v_i\sim v_{i'}}= \mathbbm1_{v_i^j\sim v_{i'}^{j}}=0$ for each $j\in [n]$. 
If $\mathbbm 1_{v_i^j\sim v_{i'}^{j'}}=1$ for some distinct $j,j' \in[n]$, then $\{v_1^j,\cdots ,v_n^j,  v_{n+1}^{j'},\cdots ,v_{2n}^{j'}\}$ forms a $K_{n,n}$, which leads to a contradiction.   Hence $\mathbbm 1_{v_i^j\sim v_{i'}^{j'}}=0$ for all distinct $j,j'\in [n]$. If $\mathbbm 1_{v_i\sim v_{i'}^j}=1$ for some $j$, then $\{v_1,\cdots ,v_n, v_{n+1}^{j},\cdots ,v_{2n}^j\}$ forms a $K_{n,n}$, a contradiction. Thus $\mathbbm 1_{v_i\sim v_{i'}^j}=0$ for each $j\in [n]$.
Therefore, all the characteristic functions are $0$. This leads to a contradiction as $\{v_i,v_i^j: i\in [n],j\in [n]\}$ forms an induced $nK_{1,n}$.
\end{proof}
Set $c_1=N_2\cdot N_3+N_2$ and $c_2=N_2$, we complete the proof.
\end{proof}

We restate \cref{finite local stable number} as follows.
\finitelocalstablenumber*
\begin{proof} 
We first prove the “only if” part. When $n = c_1 + c_2$, it holds that $(\alpha_L)_{c_1}(K_{n,n})=2n$, $(\alpha_L)_{c_1}(nK_{1,n})=n$, $(\alpha_L)_{c_1}(K_n+\overline K_n)=n$, and $(\alpha_L)_{c_1}(G_n)=n$. All of them are larger than $c_2$. This implies that they are not $\HH$-free. 
Hence  $\HH\le \{K_{n,n}, nK_{1,n}, K_n+\overline K_n, G_n\}$.

Next we prove the “if” part.
Set $$N_2=R_{2^{6n+1}}(2n),\; N_3=M\!R(N_2, 3n), \text{ and } N_1=N_2\cdot N_3+N_2,$$
where $M\!R(N_2, 3n)$ comes from Lemma \ref{bipartite Ramesy }.
\begin{claim}
Let $G$  be a graph satisfying 
\[ (\alpha_L)_{N_1}(G)\ge N_2, \]
then it contains an induced copy of $K_{n,n}, nK_{1,n}, K_n+\overline K_n$ or $ G_n$. 	
 \end{claim} 
\begin{proof}
Take $N_2$ vertices $\{v_i: i\in[N_2]\}$ satisfying $\alpha_L(v_i)\ge N_1$. Thus we can find $N_3$ stable neighbors of $v_i$ from $i=1$ to $N_2$ step by step, which are all distinct, as $N_1=N_2\cdot N_3+N_2$ is large enough. By Lemma \ref{bipartite Ramesy }, since $N_3=M\!R(N_2,3n)$, we can extract $3n$ neighbors $\{v_i^j  : j\in [3n]\}$ of $v_i$ for each $i\in [N_2]$ such that $\{v_i^j: i\in [N_2], j\in [3n]\}$ is a stable set, otherwise $K_{n,n}$ appears in $G$. Color the  auxiliary complete graph $K_{N_2}$ with $2^{6n+1}$ colors as follows. The edge $(i,i')$ with $i<i'$ is colored by 
\[(\mathbbm 1_{v_i\sim v_{i'}}, \mathbbm 1_{v_i\sim v_{i'}^1},\cdots, \mathbbm 1_{v_i\sim v_{i'}^{3n}}, \mathbbm 1_{v_{i'}\sim v_{i}^1}, \cdots, \mathbbm 1_{v_{i'}\sim v_{i}^{3n}}).\]
By Ramsey's theorem, there exists a monochromatic clique $K$ of order $2n$ in $K_{N_2}$ as $N_2=R_{2^{6n+1}}(2n)$. Without loss of generality, we may assume $V(K)=[2n]$. Consider the color of $K$. 

First suppose that  $\mathbbm 1_{v_i\sim v_{i'}}=0$. If for every $j\in [3n]$, both $\mathbbm 1_{v_i\sim v_{i'}^j}$ and $\mathbbm 1_{v_{i'}\sim v_{i}^j}: j\in [3n] $ are all $0$, then $G$ contains $nK_{1,n}$. Otherwise, by symmetry, there exists some $j\in[3n]$ such that $\mathbbm 1_{v_i\sim v_{i'}^j }=1$, then $\{v_i:i\in [n]\}\cup\{v_{i'}^j : i'\in [n+1,2n]\}$ forms a $K_{n,n}$.

Now suppose that $\mathbbm1_{v_i\sim v_{i'}}=1$. Let 
$$J=\{j\in [3n]:\mathbbm 1_{v_i\sim v_{i'}^j }=1\},\qquad J'=\{j\in [3n]:1_{v_{i'}\sim v_{i}^j }=1\}.$$ 
If $|J|\ge n$, then $\{v_i: i\in [n]\}\cup\{v_{n+1}^j: j\in J\}$ forms a $K_n+\overline K_n$. Thus $|J|\le n-1$ and $|J'|\le n-1$ by symmetry. This means that each $v_i$ has $n$  private neighbors not adjacent to any other $v_{k},\, (k\not =i)$. 
Hence $G$ contains a $G_n$ as an induced subgraph.	
\end{proof}
Set $c_1=N_1$ and $c_2=N_2$. The proof is completed.

\end{proof}
\cref{finite local component number} is restated as follows.
\finitelocalcomponentnumber*
\begin{proof}
We first prove the ``only if'' part. Set $n=c_1+c_2$. Then $(c_L)_{c_1}(K_{n,n})=2n$, $(c_L)_{c_1}(nK_{1,n})=n$ and $(c_L)_{c_1}(G_n)=n$. Thus they are not $\HH$-free, which implies that $\HH\le \{K_{n,n}, nK_{1,n}, G_n\}$.

Next we prove the ``if'' part. Set $$N_2=R_{2^{6n+1}}(2n),\; N_3=M\!R(N_2, 3n), \text{ and } N_1=N_2\cdot N_3+N_2,$$
where $M\!R(N_2, 3n)$ comes from Lemma \ref{bipartite Ramesy }.
\begin{claim}
Let $G$  be a graph satisfying 
\[ (c_L)_{N_1}(G)\ge N_2, \]
then it contains an induced copy of $K_{n,n}, nK_{1,n}$ or $ G_n$. 	
 \end{claim} 
\begin{proof}
Take $N_2$ vertices $\{v_i: i\in[N_2]\}$ satisfying $c_L(v_i)\ge N_1$. Thus we can find $N_3$ neighbors of $v_i$ lying in distinct connected components of $G[N(v_i)]$ from $i=1$ to $N_2$ step by step. These $N_2\cdot N_3$ neighbors are all distinct, as $N_1=N_2\cdot N_3+N_2$ is large enough. By Lemma \ref{bipartite Ramesy }, since $N_3=M\!R(N_2,3n)$, we can extract $3n$ neighbors $\{v_i^j  : j\in [3n]\}$ of $v_i$ for each $i\in [N_2]$ such that $\{v_i^j: i\in [N_2], j\in [3n]\}$ is a stable set, otherwise $K_{n,n}$ appears in $G$. Color the  auxiliary complete graph $K_{N_2}$ by $2^{6n+1}$ colors as follows. The edge $(i,i')$ with $i<i'$ is colored by 
\[(\mathbbm 1_{v_i\sim v_{i'}}, \mathbbm 1_{v_i\sim v_{i'}^1},\cdots, \mathbbm 1_{v_i\sim v_{i'}^{3n}}, \mathbbm 1_{v_{i'}\sim v_{i}^1}, \cdots, \mathbbm 1_{v_{i'}\sim v_{i}^{3n}}).\]
By Ramsey's theorem, there exists a monochromatic clique $K$ of order $2n$ in $K_{N_2}$ as $N_2=R_{2^{6n+1}}(2n)$. Without loss of generality, we may assume $V(K)=[2n]$. Consider the color of $K$. 

First suppose that  $\mathbbm 1_{v_i\sim v_{i'}}=0$. If for every $j\in [3n]$, both $\mathbbm 1_{v_i\sim v_{i'}^j}$ and $\mathbbm 1_{v_{i'}\sim v_{i}^j}: j\in [3n] $ are all $0$, then $G$ contains $nK_{1,n}$. Otherwise, by symmetry, there exists some $j\in[3n]$ such that $\mathbbm 1_{v_i\sim v_{i'}^j }=1$, then $\{v_i:i\in [n]\}\cup\{v_{i'}^j : i'\in [n+1,2n]\}$ forms a $K_{n,n}$.

Now suppose that $\mathbbm1_{v_i\sim v_{i'}}=1$. For fixed $i\in [2n]$, let $V_i=\big\{v_{k}:k\in [2n]-\{i\}\big\}$. Thus $V_i \subseteq N(v_i).$ 
For distinct $j,j'\in [3n]$, since $v_i^j$ cannot be connected to $v_i^{j'}$ in graph $G[N(v_i)]$, it holds that $N(v_i^j)\cap N(v_i^{j'})\cap V_i=\emptyset$. Therefore we may assume $N(v_i^j)\cap V_i=\emptyset$ for any $j\in [n]$ by the pigeonhole principle. At this time, $\{v_i,v_i^j:i\in [n],j\in [n]\}$ forms an induced $G_n$.
\end{proof}
Set $c_1=N_1$, and $c_2=N_2$, the proof is completed. 
\end{proof}

\cref{finite sdeg} is restated as follows.
\finitesdeg*
\begin{proof}
We first prove the ``only if'' part.
Set $n=c_1+c_2$. Then $(\sdeg)_{c_1}(nK_{1,n})=n$ and $(\sdeg)_{c_1}(G_n)=n$. Thus they are not $\HH$-free, which implies that $\HH\le \{ nK_{1,n}, G_n\}$.

Next we prove the ``if'' part. Set $N_2=R_{2}(n),\; N_3=N_2+n-1$, and $N_1=N_2\cdot N_3+N_2$. 
\begin{claim}
For every $\{nK_{1,n}, G_n\}$-free graph $G$, we have 
\[(\sdeg)_{N_1}(G)<N_2. \]	
\end{claim}

\begin{proof}
Suppose, for a contradiction, that there exist $N_2$ vertices $\{v_i: i\in [N_2]\}$ satisfying $\sdeg(v_i)\ge N_1=N_2\cdot N_3+N_2$. Hence there is a stable set $\{v_i^j: j\in [N_3]\} \subseteq N(v_i)$ for every $i\in[N_2]$ such that the path connecting $v_i^j$ and $v_i^{j'}$ must contain $v_i$. Moreover, all vertices $v_i,v_i^j$ for $i\in[N_2],j\in [N_3]$ are distinct.
If there are two vertices $v_i^j$ and $v_i^{j'}$ satisfying 
\[N(v_i^j)\cap \{v_{i'}, v_{i'}^k: k\in [N_3]\}\not =\emptyset \text{ and } N(v_{i}^{j'})\cap \{v_{i'}, v_{i'}^k: k\in [N_3]\}\not =\emptyset ,\]
for some $i'\neq i$, then $v_i^j$ can be connected to $v_i^{j'}$ through the vertex set $\{v_{i'},v_{i'}^k: k\in [N_3]\}$, a contradiction. Therefore, by retaining only $N_3-(N_2-1)=n$ vertices $v_i^j$ satisfying $N(v_i^j)\cap \{v_{i'}, v_{i'}^k: i'\in [N_2]\setminus\{i\}, k\in [N_3]\} =\emptyset$ from $\{v_i^j:j\in [N_3]\}$ for any $i\in [N_2]$,
we may assume 
\[N(v_i^j)\cap \{v_{i'}, v_{i'}^k: k\in [n]\} =\emptyset,\; \text{for all $i\not = i'$ and all } j\in [n].\] 
Note that the set $V=\{v_i: i\in [N_2]\}$ has order $N_2=R_2(n)$. By Ramsey's Theorem again, if $V$ has a clique of order $n$, we have a $G_n$ in $G$; otherwise, $V$ has a stable set of order $n$, and we have an $nK_{1,n}$ in $G$. 
\end{proof}

Set $c_1=N_1$, and $c_2=N_2$, the proof is completed. 
\end{proof}

%Note that the proof previously appeared can be adjusted slightly and the bound can be made smaller. For example, we can use Theorem \ref{domination number} to handle theorems in section $2$ as the same way as Theorem \ref{local connected}.
%We can  require that the graphs in $\HH$ use different positive integer, such as $\HH\le \{K_{1,n_1}, K_{n_2}, P_{n_3}\}$ in Theorem \ref{cut vertex}. We can also require that the constant $c_1(\HH)$ and $c_2(\HH)$ be the same in Theorem \ref{finite deg} since we can replace  $c_1(\HH)$ and $c_2(\HH)$ by the same constant $c_1(\HH)+c_2(\HH)$.

 %However, it has no meaning except making proofs seem more complicated. As for Theorem \ref{finite deg}, Theorem \ref{finite local stable number} and Theorem \ref{finite adh}, we don't require $G$ is connected, because there is no typical way to connect the $n$ discrete $K_{1,n}$. For a tree with many leaves, replace the leaves by $K_{1,n}$, then we get a counterexample. However, we have no good way to forbid this case. The more explanation, see \cite{c20}, \cite{l17}.\bigskip

%\noindent{\bf Acknowledgment:} {The work was supported by the National Natural Science Foundation of China (No. 12071453), the National Key R and D Program of China(2020YFA0713100), the Anhui Initiative in Quantum Information Technologies (AHY150200)  and the Innovation Program for Quantum Science and Technology, China (2021ZD0302904).}
%\noindent{\bf Data Availability:} Data sharing is not applicable to this article as no datasets were generated or analyzed during the current study.

\medskip
\noindent{\bf Acknowledgements:} This work was supported by the National Key R\&D Program of China
		(2023YFA1010203), the National Natural Science Foundation of China
		(12471336), and the  Quantum Science and Technology--National Science and Technology Major Project (2021ZD0302902).


\begin{thebibliography}{99}
\bibitem{ado}
S. Allred, G. Ding, and B. Oporowski, Unavoidable induced subgraphs of large 2-connected graphs, SIAM J. Discrete Math. {\bf 37} (2023), no.2, 684--698.

%\bibitem{al}B. Alecu, V. Alekseev, A. Atminas, V. Lozin, and V. Zamaraev, Graph parameters, implicit representation and factorial properties of graphs, arXiv:2303.04453.

\bibitem{alr}
A. Atminas, V. V. Lozin, and I. Razgon, Linear time algorithm for computing a small biclique in graphs without long induced paths, In Algorithm theory-SWAT 2012, volume 7357 of Lecture Notes in Comput. Sci. Springer, Heidelberg, 142--152.

%\bibitem{cf}S. Chiba and M. Furuya, Ramsey-type results for path covers and path partitions, The Electronic Journal of Combinatorics 29(4), (2022).

%\bibitem{cfkp}I. Choi, M. Furuya, R. Kim and B. Park, A Ramsey-type theorem for the matching number regarding connected graphs, Discrete Math. 343 (2020).

\bibitem{ckos} M. Chudnovsky, R. Kim, S. Oum, and P. Seymour, Unavoidable induced subgraphs in large graphs with no homogeneous sets, J. Combin. Theory Ser. B {\bf 118} (2016), 1--12.

\bibitem{d}
R. Diestel, Graph theory, GTM 173, Springer, Berlin, fifth edition, (2017).

\bibitem{dc}
G. Ding and P. Chen, Unavoidable doubly connected large graphs, Discrete Mathematics {\bf 280} (2004), 1--12.

\bibitem{es}
D. Eppstein and E.S. Spiro, The $h$-index of a graph and its application to dynamic sub-graph statistics, J. Graph Algorithms and Applications {\bf 16}, (2012) 543-567.

\bibitem{e}
P. Erd\H{o}s, Graph theory and probability, Canad. J. Math. {\bf 11} (1959), 34--38.

\bibitem{f}
M. Furuya, Forbidden subgraphs for constant domination number, Discrete Math. Theor. 20:1, (2018).

\bibitem {g}
A. Gy\'arf\'as. Problems from the world surrounding perfect graphs. In Proc. Int. Conf. on Comb. Analysis and Applications (Pokrzywna, 1985). Zastos. Mat. {\bf 19}, (1987) 413--441. 

\bibitem{grs}
F. Galvin, I. Rival, and B. Sands, A Ramsey-type theorem for traceable graphs. J. Combin. Theory, Ser. B, {\bf 33} (1982) 7-16. 

\bibitem{hhbgm}
M. S. Handcock, D. Hunter, C. T. Butts, S. M. Goodreau, and M. Morris, Statnet: An R package for the Statistical Modeling of Social Networks. Web page http://www.csde.washington.edu/statnet, (2003).

\bibitem{kp}
H. A. Kierstead and S. G. Penrice, Radius two trees specify $\chi$-bounded classes. J. Graph Theory {\bf 18} (1994) 119--129.

\bibitem {lr}
V. Lozin and I. Razgon. Tree-width dichotomy. European J. Combin., {\bf 103} (2022), 103517.

\bibitem {r}
F. P. Ramsey, On a Problem of Formal Logic, Proc. London Math. Soc. (2), 30(4), (1929) 264--286.

\bibitem{sh}J. Sun and X. Hou, A Ramsey-type theorem on deficiency, J. Graph Theory, {\bf 110} (2025), 313--321.

\bibitem{sss}
A. Scott, P. Seymour, and S. Spirkl, Polynomial bounds for chromatic number. I. Excluding a biclique and an induced tree, J. Graph Theory {\bf 102} (2023) 458--471.

\bibitem{s81} D. P. Sumner, Subtrees of a graph and chromatic number. in The Theory
and Applications of Graphs, John Wiley $\&$ Sons, New York (1981), 557--576.


\end{thebibliography}
\end{document}